\documentclass{ifacconf}

\usepackage{epsfig} 
\usepackage{amsmath} 
\usepackage{amssymb}  
\usepackage{url}

\usepackage{enumerate}

\begin{document}

\begin{frontmatter}

\title{Improving Fast Dual Ascent for MPC - Part I: The Distributed Case\thanksref{footnoteinfo}}

\thanks[footnoteinfo]{
During the preparation of this paper, the author was a
member of the LCCC Linnaeus Center at Lund University. 
Financial support from the Swedish Research Council for the
author's Postdoctoral studies at Stanford University is gratefully
acknowledged. Eric Chu is also gratefully acknowledged for
constructive feedback.}

\author[Pontus]{Pontus Giselsson} 

%
\address[Pontus]{Electrical Engineering, Stanford University\\
   (e-mail: pontusg@stanford.edu).}                                              



\begin{abstract}

In dual decomposition, the dual to an optimization problem with a
specific structure is solved in distributed fashion using
(sub)gradient and recently also fast gradient methods. The traditional
dual decomposition suffers from two main short-comings. The first is
that the convergence is often slow, although fast gradient
methods have significantly improved the situation. The
second is that computation of the
optimal step-size requires centralized computations, which hinders a fully
distributed 
implementation of the algorithm. In this paper, the first issue is addressed
by providing a tighter characterization of the dual function than what
has previously been reported in the literature. Then a distributed and
a parallel
algorithm are presented in which the provided dual function approximation is
minimized in each step. Since the approximation is more accurate than
the approximation used in standard and fast dual decomposition, the
convergence properties are improved. For the second issue, we
extend a recent result to allow for a fully distributed parameter
selection in the algorithm. Further, we show how to apply the
proposed algorithms to optimization problems arising in distributed
model predictive control (DMPC) and show that the proposed distributed algorithm enjoys
distributed reconfiguration, i.e. \emph{plug-and-play}, in the DMPC context.

\end{abstract}

\end{frontmatter}


\section{Introduction}

Optimization problems with a separable cost and sparse constraints can
be solved in distributed fashion by distributed optimization
algorithms. Some distributed algorithms
exploit the property that the (sub)gradient to the dual of
such optimization problems can
be computed in distributed fashion, which enables for
distributed implementation of dual (sub)gradient algorithms. This approach is
referred to as dual decomposition and
originates from \cite{everett,Danzig+61,Benders}. The use of
sub-gradient or gradient methods to 
solve the dual problem usually results in poor convergence properties of
the algorithm. As a remedy to this, a dual Newton method was presented
in \cite{dualNewtonCG} where the dual problem is solved in distributed
fashion using a Newton method. The Newton step is computed in
distributed fashion using a distributed implementation of a conjugate
gradient method. Another distributed algorithm was presented in \cite{BlockSplitBoyd},
which is based on the alternating direction method of multipliers
(ADMM, see \cite{BoydDistributed}), and solves a more general class of problems
than dual decomposition and the dual Newton method in \cite{dualNewtonCG}.
In \cite{gisAutomatica}, another recent attempt to improve
the convergence of distributed 
algorithms was presented. It relies on using fast gradient methods in dual
decomposition. These fast gradient methods were originally
presented in \cite{Nesterov1983} in the early 80's. These methods
rendered no or little attention the following decades but became
increasingly studied from the mid 00's. Since then, the fast
gradient method has been extended and
generalized in several directions, see e.g.
\cite{BecTab_FISTA:2009,NesterovLectures,Tseng_acc:2008,Nesterov2005}. The main
benefit of fast gradient methods is that, with negligible
increase in iteration complexity, the convergence rate is improved
from $O(1/k)$ for standard gradient methods to $O(1/k^2)$, where $k$ is the iteration number.
Obviously, the use of fast gradient methods in dual decomposition
instead of standard gradient methods has considerably improved
the convergence properties. However,
in many applications further improvements are
necessary for realistic implementation. In this paper, we propose
dual decomposition like algorithms that have further improved
convergence properties.

In a general form, fast gradient methods can be applied to problems
consisting of a sum of two functions.
The prerequisites for these functions are that one is convex and differentiable and has a
Lipschitz continuous gradient, while the other is proper, closed, and
convex. The former properties are equivalent to the
existence of a
quadratic upper bound with the same curvature in all directions (defined by the
Lipschitz constant) to the function. In gradient and
fast gradient
methods, this quadratic upper bound is used as an approximation to the
function. This approximation plus the closed, proper, convex function
is minimized in every step of the algorithm. If
the quadratic upper bound does not well
approximate the function, slow convergence
properties are expected. By instead letting the quadratic upper bound
have different curvature in different directions, a closer fit between
the bound and the function can be obtained. For an appropriate choice
of non-uniform quadratic upper bound, this can significantly improve
the convergence properties of fast gradient methods. 
The key result of this paper is a characterization of the set of
matrices that can be used to 
describe a quadratic upper bound to the convex negative dual function, in the case of
strongly convex primal cost function. This result generalizes previous
results, e.g. \cite{Nesterov2005}, where a Lipschitz
constant to the dual gradient is quantified. As a consequence of the
presented result, quadratic upper bounds with different curvature in different
directions can be used in dual decomposition methods.

In this paper, we propose two improved dual decomposition algorithms
based on the previously mentioned key result, one parallel and one
distributed version. For both algorithms, the
matrix that describes the quadratic upper bound to the dual function
must be chosen.
In the parallel version, there are no restrictions on the structure of
that matrix, while in the distributed algorithm, the matrix must be
block diagonal to facilitate a distributed implementation. In 
fast dual decomposition,
this matrix is traditionally chosen as
the reciprocal of the Lipschitz constant to the dual gradient times the
identity matrix. By allowing for more flexibility in the
matrix structure, the
shape of the minimized function can be better captured. In this paper,
we also show how to
compute a matrix that, when used as basis for the quadratic upper bound in the dual
decomposition algorithm, can significantly improve the convergence. 

Besides convergence issues in dual decomposition, there is the
issue of computing the step-size. The optimal choice requires the
computation of the 2-norm of a system-wide matrix. This cannot
straight-forwardly be done
using distributed computations. However, approximations to this norm
can be computed in distributed fashion with centralized coordination.
In this paper, we
extend a recent result in \cite{BeckEtal2013} to enable a fully distributed
initialization procedure for our distributed algorithm. The
initialization selects a
block diagonal matrix that describes the
quadratic upper bound using local computations and communication only.
For our parallel algorithm, the initialization need not be distributed
since the algorithm needs all data to be accessible in a centralized
unit.

In distributed model predictive control (DMPC), dual decomposition
techniques have been used to distribute the computations over the
subsystems
\cite{NegenbornPhd,doan2011,gisAutomatica}.
Although the use of fast gradient methods in dual decomposition have
significantly improved the
convergence, see \cite{gisAutomatica}, it is not enough for
realistic implementation in a distributed control system. In
\cite{gisACC2013a}, a generalized version of dual decomposition was
presented that allows for different curvature in different directions
in the quadratic upper bound that is minimized in every iteration of
the algorithm. This gives a significantly reduced
number of iterations. The algorithm in \cite{gisACC2013a} is restricted
to problems having a quadratic cost, linear equality constraints, and
linear inequality constraints. Dual variables for all these constraints are
introduced, which results in the dual problem being a quadratic program. The algorithm in
this paper is an extension and 
generalization of the algorithm in \cite{gisACC2013a} that allows for
any (local) convex inequality constraints. Also, only the equality
constraints are dualized in this paper. These changes give rise to completely
different technicalities since the dual function is implicitly
defined though an optimization problem.

A feature of
DMPC is that similar optimization problems are repeatedly solved
online. This implies that much offline computational
effort can be devoted to parameter selection in the algorithm to
improve the online convergence. In this paper, the offline computational effort is
devoted to choose a matrix that describes the quadratic upper bound
to the negative dual function. The numerical evaluation suggests that
this can significantly reduce the number of iterations in the
algorithm compared to dual decomposition using fast gradient methods,
and compared to the dual Newton method in \cite{dualNewtonCG}.
Besides favorable convergence
properties, the presented distributed algorithm enjoys distributed
configuration and reconfiguration, commonly referred to as
plug-and-play. Distributed reconfiguration or plug-and-play
is the property that
if a subsystem is added to (or removed from) the system, only
neighboring subsystems need to be invoked to reconfigure the
algorithm for the new setup.

This paper is an extension of \cite{gis2014IFACdist}, and is the first
paper in a series of two on improving fast dual ascent for model
predictive control, where \cite{gis2014AutPart2} is the second.

\section{Preliminaries and Notation}

\subsection{Notation}

We denote by $\mathbb{R}$, $\mathbb{R}^n$, $\mathbb{R}^{m\times n}$,
the sets of real numbers, vectors, and matrices.
$\mathbb{S}^n\subseteq\mathbb{R}^{n\times n}$ is the set of symmetric matrices, and 
$\mathbb{S}_{++}^n\subseteq\mathbb{S}^n$,
$[\mathbb{S}_{+}^n]\subseteq\mathbb{S}^n$, are the sets of positive [semi] 
definite matrices. Further, $L\succeq M$ and $L\succ M$
where $L,M\in\mathbb{S}^n$ denotes
$L-M\in\mathbb{S}_{+}^n$ and $L-M\in\mathbb{S}_{++}^n$ respectively.
We also use notation
$\langle x,y\rangle=x^Ty$, $\langle
x,y\rangle_H=x^THy$, $\|x\|_2=\sqrt{x^Tx}$, and $\|x\|_H =
\sqrt{x^THx}$. Finally, $I_{\mathcal{X}}$ denotes the indicator
function for the set $\mathcal{X}$, i.e.
$I_{\mathcal{X}}(x)\triangleq\left\{\begin{smallmatrix}0,~&x\in\mathcal{X}\\
\infty, &{\rm{else~}}\end{smallmatrix}\right.$.

\subsection{Preliminaries}

In this section, we introduce generalizations of already well used
concepts. We generalize the notion of strong convexity as well as
the notion of Lipschitz continuity of the gradient of convex
functions. We also define conjugate
functions and state a known result on dual properties of a
function and its conjugate.

For differentiable and convex functions
$f~:~\mathbb{R}^n\to\mathbb{R}$ that have a Lipschitz continuous
gradient with constant $L$, we have that
\begin{equation}
\|\nabla f(x_1)-\nabla f(x_2)\|_2\leq L\|x_1-x_2\|_2
\label{eq:gradLipschitz}
\end{equation}
holds for all $x_1,x_2\in\mathbb{R}^n$.
This is equivalent to that
\begin{equation}
f(x_1)\leq f(x_2)+\langle \nabla
f(x_2),x_1-x_2\rangle+\frac{L}{2}\|x_1-x_2\|_2^2
\label{eq:standardQuadBound}
\end{equation}
holds for all
$x_1,x_2\in\mathbb{R}^n$ \cite[Theorem 2.1.5]{NesterovLectures}.
In this paper, we allow for a generalized version of the quadratic
upper bound \eqref{eq:standardQuadBound} to $f$, namely that
\begin{equation}
f(x_1)\leq f(x_2)+\langle \nabla
f(x_2),x_1-x_2\rangle+\frac{1}{2}\|x_1-x_2\|_{\mathbf{L}}^2
\label{eq:quadBound}
\end{equation}
holds for all $x_1,x_2\in\mathbb{R}^n$ where
$\mathbf{L}\in\mathbb{S}_+^{n}$. The
bound \eqref{eq:standardQuadBound} is obtained by setting
$\mathbf{L}=LI$ in \eqref{eq:quadBound}. 
\begin{rem}
For concave functions $f$, i.e. where $-f$ is convex, the Lipschitz
condition \eqref{eq:gradLipschitz} is equivalent to that the following
quadratic lower bound
\begin{equation}
f(x_1)\geq f(x_2)+\langle \nabla
f(x_2),x_1-x_2\rangle-\frac{L}{2}\|x_1-x_2\|_2^2
\label{eq:standardQuadBoundConcave}
\end{equation}
holds for all $x_1,x_2\in\mathbb{R}^n$. The generalized counterpart
naturally becomes that 
\begin{equation}
f(x_1)\geq f(x_2)+\langle \nabla
f(x_2),x_1-x_2\rangle-\frac{1}{2}\|x_1-x_2\|_{\mathbf{L}}^2
\label{eq:quadBoundConcave}
\end{equation}
holds for all $x_1,x_2\in\mathbb{R}^n$.
\end{rem}
Next, we state a
Lemma on equivalent characterizations of the condition \eqref{eq:quadBound}.
\begin{lem}
Assume that $f~:~\mathbb{R}^n\to\mathbb{R}$
is convex and differentiable. The condition that 
\begin{equation}
f(x_1)\leq f(x_2)+\langle
\nabla f(x_2),x_1-x_2\rangle+\frac{1}{2}\|x_1-x_2\|_{\mathbf{L}}^2
\label{eq:quadBoundLem}
\end{equation}
holds for some
$\mathbf{L}\in\mathbb{S}_{+}^n$ and all $x_1,x_2\in\mathbb{R}^n$
is equivalent to that 
\begin{equation}
\langle\nabla f(x_1)-\nabla f(x_2),x_1-x_2\rangle \leq \|x_1-x_2\|_{\mathbf{L}}^2.
\label{eq:quadBoundCondition}
\end{equation}
holds for all $x_1,x_2\in\mathbb{R}^n$.
\label{lem:quadUBequiv}
\end{lem}
\begin{pf}
To show the equivalence, we introduce the function
$g(x):=\frac{1}{2}x^T\mathbf{L}x-f(x)$. According to \cite[Theorem
2.1.3]{NesterovLectures} and since $g$ is differentiable,
$g~:~\mathbb{R}^n\to\mathbb{R}$ is convex if and only if $\nabla g$ is
monotone. The function $g$ is convex if and only if 
\begin{align*}
&g(x_1)\geq g(x_2)+\langle\nabla g(x_2),x_1-x_2\rangle=\\
&=\frac{1}{2}x_2^T\mathbf{L}x_2-f(x_2)+\langle \mathbf{L}x_2-\nabla
f(x_2),x_1-x_2\rangle\\
&=-f(x_2)-\langle \nabla
f(x_2),x_1-x_2\rangle-\tfrac{1}{2}\|x_1-x_2\|_\mathbf{L}^2+\tfrac{1}{2}x_1^T\mathbf{L}x_1.
\end{align*}
Noting that $g(x_1) = \frac{1}{2}x_1^T\mathbf{L}x_1-f(x_1)$ gives the
negated version of \eqref{eq:quadBoundLem}.

Monotonicity of $\nabla g$ is equivalent to 
\begin{align*}
0&\leq \langle \nabla g(x_1)-\nabla g(x_2),x_1-x_2\rangle\\
&=\langle \mathbf{L}x_1-\nabla f(x_1)-\mathbf{L}x_2+\nabla f(x_2),x_1-x_2\rangle\\
&=\|x_1-x_2\|_{\mathbf{L}}^2-\langle \nabla f(x_1)-\nabla f(x_2),x_1-x_2\rangle.
\end{align*}
Rearranging the terms gives \eqref{eq:quadBoundCondition}. This
concludes the proof.
\end{pf}

Next, we state the corresponding result for concave functions.
\begin{cor}
Assume that $f~:~\mathbb{R}^n\to\mathbb{R}$
is concave and differentiable. The condition that 
\begin{equation}
f(x_1)\geq f(x_2)+\langle
\nabla f(x_2),x_1-x_2\rangle-\frac{1}{2}\|x_1-x_2\|_{\mathbf{L}}^2
\label{eq:quadBoundCorConcave}
\end{equation}
holds for some
$\mathbf{L}\in\mathbb{S}_{+}^n$ and all $x_1,x_2\in\mathbb{R}^n$
is equivalent to that 
\begin{equation}
\langle\nabla f(x_1)-\nabla f(x_2),x_2-x_1\rangle \leq \|x_1-x_2\|_{\mathbf{L}}^2.
\label{eq:quadBoundConditionConcave}
\end{equation}
holds for all $x_1,x_2\in\mathbb{R}^n$.
\label{cor:quadUBequivConcave}
\end{cor}
\begin{pf}
The proof follows directly from $-f$ being convex and applying
Lemma~\ref{lem:quadUBequiv}. 
\end{pf}

The standard definition of a differentiable and strongly convex 
function $f~:~\mathbb{R}^n\to\mathbb{R}$ is that it satisfies
\begin{equation}
f(x_1)\geq f(x_2)+\langle \nabla
f(x_2),x_1-x_2\rangle+\frac{\sigma}{2}\|x_1-x_2\|_2^2 
\label{eq:strConvStandardDef}
\end{equation}
for any $x_1,x_2\in\mathbb{R}^n$, where the modulus $\sigma\in\mathbb{R}_{++}$
describes a lower bound of the curvature of the function. In this
paper, the definition \eqref{eq:strConvStandardDef} is generalized to
allow for a quadratic lower bound with different curvature in different
directions.
\begin{defn}
A differentiable function $f~:~\mathbb{R}^n\to\mathbb{R}$
is \emph{strongly convex with matrix H} if and only if
\begin{equation*}
f(x_1)\geq f(x_2)+\langle \nabla
f(x_2),x_1-x_2\rangle+\frac{1}{2}\|x_1-x_2\|_H^2 
\end{equation*}
holds for all $x_1,x_2\in\mathbb{R}^n$, where
$H\in\mathbb{S}_{++}^{n}$.
\label{def:strConv}
\end{defn}
\begin{rem}
The traditional definition of strong convexity
\eqref{eq:strConvStandardDef} is obtained from
Definition~\ref{def:strConv} by setting $H=\sigma I$.
\end{rem}
\begin{lem}
Assume that $f~:~\mathbb{R}^n\to\mathbb{R}$
is differentiable and
strongly convex with matrix $H$. The condition that
\begin{equation}
f(x_1)\geq f(x_2)+\langle
\nabla f(x_2),x_1-x_2\rangle+\frac{1}{2}\|x_1-x_2\|_{H}^2
\label{eq:strConv}
\end{equation}
holds for all $x_1,x_2\in\mathbb{R}^n$
is equivalent to that
\begin{equation}
\langle\nabla f(x_1)-\nabla f(x_2),x_1-x_2\rangle \geq \|x_1-x_2\|_{H}^2
\label{eq:strConvCondition}
\end{equation}
holds for all $x_1,x_2\in\mathbb{R}^n$.
\label{lem:strConvCondition}
\end{lem}
\begin{pf}
To show the equivalence, we introduce the function
$g(x):=f(x)-\frac{1}{2}x^THx$ and proceed similarly to in the proof
of Lemma \eqref{lem:quadUBequiv}. According to \cite[Theorem
2.1.3]{NesterovLectures} and since $g$ is differentiable,
$g~:~\mathbb{R}^n\to\mathbb{R}$ is convex if and only if $\nabla g$ is
monotone. The function $g$ is convex if and only if 
\begin{align*}
&g(x_1)\geq g(x_2)+\langle\nabla g(x_2),x_1-x_2\rangle=\\
&=f(x_2)-\frac{1}{2}x_2^THx_2+\langle \nabla
f(x_2)-Hx_2,x_1-x_2\rangle\\
&=f(x_2)+\langle \nabla
f(x_2),x_1-x_2\rangle+\tfrac{1}{2}\|x_1-x_2\|_H^2-\tfrac{1}{2}x_1^THx_1.
\end{align*}
Noting that $g(x_1) = f(x_1)-\frac{1}{2}x_1^THx_1$ gives \eqref{eq:strConv}.

Monotonicity of $\nabla g$ is equivalent to 
\begin{align*}
0&\leq \langle \nabla g(x_1)-\nabla g(x_2),x_1-x_2\rangle\\
&=\langle \nabla f(x_1)-Hx_1-\nabla f(x_2)+\mathbf{L}x_2,x_1-x_2\rangle\\
&=\langle \nabla f(x_1)-\nabla f(x_2),x_1-x_2\rangle-\|x_1-x_2\|_{H}^2.
\end{align*}
Rearranging the terms gives \eqref{eq:strConvCondition}. This
concludes the proof.
\end{pf}

The condition \eqref{eq:strConv} is a quadratic lower bound on the function value,
while the condition \eqref{eq:quadBound} is a quadratic upper bound on the function
value. These two properties are linked through the conjugate function
\begin{align*}
f^{\star}(y) &\triangleq \sup_x \left\{y^Tx-f(x)\right\}.
\end{align*}
More precisely, we have the following result.
\begin{prop}
Assume that $f~:~\mathbb{R}^n\to\mathbb{R}\cup\{\infty\}$ is closed,
proper, and strongly convex with
modulus $\sigma$ on the relative interior of its domain. Then the
conjugate function
$f^{\star}$ is convex and differentiable, and
$\nabla f^{\star}(y) = x^\star(y)$, where $x^\star(y)=\arg\max_x 
\left\{y^Tx-f(x)\right\}$. Further, $\nabla f^{\star}$ is Lipschitz continuous with 
constant $L = \frac{1}{\sigma}$.
\label{prp:strConvSmooth}
\end{prop}
A straight-forward generalization is given by the chain-rule and was
proven in \cite[Theorem 1]{Nesterov2005} (which also proves the less
general Proposition~\ref{prp:strConvSmooth}).
\begin{cor}
Assume that $f~:~\mathbb{R}^n\to\mathbb{R}\cup\{\infty\}$ is closed,
proper, and strongly convex with
modulus $\sigma$ on the relative interior of its domain. Further, define
$g^\star(y) \triangleq f^\star(Ay)$. Then
$g^\star$ is convex and differentiable, and $\nabla g^\star(y) =
A^Tx^\star(Ay)$, where $x^\star(Ay)=\arg\max_x 
\left\{(Ay)^Tx-f(x)\right\}$. Further, $\nabla g^{\star}$
is Lipschitz continuous with 
constant $L = \frac{\|A\|_2^2}{\sigma}$.
\label{cor:strConvSmooth}
\end{cor}
For the case when $f(x) = \frac{1}{2}x^THx+g^Tx$, i.e. $f$ is a
quadratic, a tighter Lipschitz constant to $\nabla g^\star(y)=\nabla f^\star(Ay)$ was provided
in \cite[Theorem 7]{RichterMMOR2013}, namely
$L=\|AH^{-1}A^T\|_2$.

\section{Problem formulation}

We consider optimization problems of the form
\begin{equation}
\begin{tabular}[t]{ll}
minimize & $f(x)+h(x)+g(Bx)$\\
subject to & $Ax=b$
\end{tabular}
\label{eq:primProb}
\end{equation}
where the decision variables are partitioned as
$x=(x_1,\ldots,x_M)\in\mathbb{R}^{n}$ where
$x_i\in\mathbb{R}^{n_i}$, the cost functions are
separable, i.e., $f(x) = \sum_{i=1}^M
f_i(x_i)$, $h(x) = \sum_{i=1}^M h_i(x_i)$, and $g(Bx) = \sum_{i=1}^M
g_i(B_ix)$, where $B\in\mathbb{R}^{p\times n}$ and
$B_i\in\mathbb{R}^{p_i\times n}$ for $i=\{1,\ldots,M\}$, are partitioned as
\begin{align*}
B &= \begin{bmatrix}
B_{1} \\
\vdots \\
B_{M} 
\end{bmatrix},&
B_i &= \begin{bmatrix}
B_{i1} & \cdots & B_{iM}
\end{bmatrix}
\end{align*}
where $B_{ij}\in\mathbb{R}^{p_i\times n_j}$ for all $i\in\{1,\ldots,M\}$
and $j\in\{1,\ldots,M\}$.
Further, $A\in\mathbb{R}^{m\times n}$ and $b\in\mathbb{R}^m$, are partitioned as
\begin{align*}
A &= \begin{bmatrix}
A_{11} & \cdots & A_{1M}\\
\vdots & \ddots & \vdots\\
A_{M1} & \cdots & A_{MM}
\end{bmatrix},&
b &= \begin{bmatrix}
b_{1} \\
\vdots \\
b_{M} 
\end{bmatrix}
\end{align*}
where $A_{ij}\in\mathbb{R}^{m_i\times n_j}$ for all $i\in\{1,\ldots,M\}$
and $j\in\{1,\ldots,M\}$ and $b_i\in\mathbb{R}^{m_i}$ for all
$i\in\{1,\ldots,M\}$.
We assume for all $i\in\{1,\ldots,M\}$, that $A_{ij}=0$ and $B_{ij}=0$
for
some $j\in\{1,\ldots,M\}$, i.e., that the $A$ and $B$ matrices are block sparse.
The sparsity structure induced by this
assumption is represented by the sets $\mathcal{N}_i$ and
$\mathcal{M}_i$, where $\mathcal{N}_i$ contains indices for non-zero
blocks of block row $i$ and $\mathcal{M}_i$ contains indices for non-zero
blocks of block column $i$. More precisely, we have
\begin{align*}
\mathcal{N}_i&=\left\{j\in\{1,\ldots,M\}~|~A_{ij}
  \neq 0 {\hbox{ and }} B_{ij} \neq 0\right\},\\
\mathcal{M}_i&=\left\{j\in\{1,\ldots,M\}~|~A_{ji}
  \neq 0 {\hbox{ and }} B_{ji} \neq 0\right\}.
\end{align*}
We also introduce concatenated matrices
$A_{\mathcal{N}_i}\in\mathbb{R}^{m_i\times n_{\mathcal{N}_i}}$, where
$n_{\mathcal{N}_i} = \sum_{j\in\mathcal{N}_i} n_j$, that 
contain all non-zero sub-matrices $A_{ij}$, e.g., if $\mathcal{N}_1
= \{1,2,6\}$
then $A_{\mathcal{N}_1}=[A_{11}~A_{12}~A_{16}]$. Similarly, we
introduce $A_{\mathcal{M}_i}\in\mathbb{R}^{m_{\mathcal{M}_i}\times
  n_i}$, where
$m_{\mathcal{M}_i} = \sum_{j\in\mathcal{M}_i} m_j$; if $\mathcal{M}_1= \{1,4,6\}$, then
$A_{\mathcal{M}_1}=[A_{11}^T~A_{41}^T~A_{61}^T]^T$. This notation
is used for all matrices that have a block structure as
specified by $\mathcal{N}_i$ and $\mathcal{M}_i$, e.g.,
$B_{\mathcal{N}_i}\in\mathbb{R}^{p_i\times n_{\mathcal{N}_i}}$ and
$B_{\mathcal{M}_i}\in\mathbb{R}^{p_{\mathcal{M}_i}\times
  n_i}$ where $p_{\mathcal{M}_i} = \sum_{j\in\mathcal{M}_i} p_j$, are defined equivalently.
We also introduce
consistent notation for the variables, namely
$x_{\mathcal{N}_i}\in\mathbb{R}^{n_{\mathcal{N}_i}}$, i.e.
$x_{\mathcal{N}_1} = (x_1,x_2,x_6)$ in the above example. This
implies that $\sum_{j\in\mathcal{N}_i} A_{ij}x_j =
A_{\mathcal{N}_i}x_{\mathcal{N}_i}$ and $B_{i}x =
B_{\mathcal{N}_i}x_{\mathcal{N}_i}$.

\begin{rem}
Note that some sub-matrices of $A_{\mathcal{N}_i}$,
$A_{\mathcal{M}_i}$, $B_{\mathcal{N}_i}$, and $B_{\mathcal{M}_i}$ may
be zero due to the construction of $\mathcal{N}_i$ and
$\mathcal{M}_i$. We allow this for notational convenience.
\end{rem}

The preceding assumptions and the introduced notation imply that the
optimization problem \eqref{eq:primProb} can
equivalently be written
\begin{equation}
\begin{tabular}[t]{lll}
minimize & \multicolumn{2}{l}{$\displaystyle \sum_{i=1}^M \left\{f_i(x_i)+h_i(x_i)+g_i(y_i)\right\}$}\\
subject to & $A_{\mathcal{N}_i}x_{\mathcal{N}_i}=b_i$, &
$i=\{1,\ldots,M\}$\\
& $B_{\mathcal{N}_i}x_{\mathcal{N}_i} = y_i$, &
$i=\{1,\ldots,M\}$
\end{tabular}
\label{eq:primProbSep}
\end{equation}
Throughout this paper we assume the following.
\begin{assum}~
\begin{enumerate}[(a)]
\item The functions $f_i~:~\mathbb{R}^{n_i}\to\mathbb{R}$ are strongly
  convex with matrix $H_i\in\mathbb{S}_{++}^{n_i}$.
\item The extended valued functions
  $h_i~:~\mathbb{R}^{n_i}\to\mathbb{R}\cup\{\infty\}$ and
  $g_i~:~\mathbb{R}^{p_i}\to\mathbb{R}\cup\{\infty\}$ are proper,
  closed, and convex.
\item The matrix $A\in\mathbb{R}^{m\times n}$ has full row rank.
\end{enumerate}
\label{ass:probAss}
\end{assum}
\begin{rem}
Assumption~\ref{ass:probAss}(a) implies that $f = \sum_{i=1}^M f_i$
is strongly convex with matrix $H$, where 
\begin{equation}
H :=
{\rm{blkdiag}}(H_1,\ldots,H_M).
\label{eq:Hdef}
\end{equation}
Assumption~\ref{ass:probAss}(b) is satisfied if, e.g., $h_i$ and $g_i$ are
indicator functions to convex constraint sets. If Assumption~\ref{ass:probAss}(c) is
not satisfied, redundant equality constraints can, without affecting
the solution of \eqref{eq:primProb}, be removed to satisfy the assumption.
\end{rem}

To form the dual problem, we introduce dual
variables $\lambda=(\lambda_1,\ldots,\lambda_M)\in\mathbb{R}^m$
where $\lambda_i\in\mathbb{R}^{m_i}$, and $\mu =
(\mu_1,\ldots,\mu_M)\in\mathbb{R}^p$ where
$\mu_i\in\mathbb{R}^{p_i}$. We also introduce a notation for
dual variables that correspond to the concatenated matrices
$A_{\mathcal{M}_i}$ and $B_{\mathcal{M}_i}$, namely
$\lambda_{\mathcal{M}_i}\in\mathbb{R}^{m_{\mathcal{M}_i}}$ and
$\mu_{\mathcal{M}_i}\in\mathbb{R}^{m_{\mathcal{M}_i}}$
respectively. In the
above example with $\mathcal{M}_1 = \{1,4,6\}$ we get 
$\lambda_{\mathcal{M}_1} = (\lambda_1,\lambda_4,\lambda_6)$
and $\mu_{\mathcal{M}_1} = (\mu_1,\mu_4,\mu_6)$. 
This gives the following Lagrange dual problem
\begin{align}
\nonumber & \displaystyle\sup_{\lambda,\mu}\inf_{x,y}\left\{
\displaystyle f(x)+h(x)+\lambda^T(Ax-b)+g(y)+\mu^T(Bx-y)\right\}\\
\nonumber &=\displaystyle\sup_{\lambda,\mu}\inf_{x,y} \sum_{i=1}^M\Big\{
f_i(x_i)+h_i(x_i)+\lambda_i^T\left(A_{\mathcal{N}_i}x_{\mathcal{N}_i}-b_i\right) \\
\nonumber&\qquad\qquad\qquad\qquad\qquad\quad+g_i(y_i)+\mu_i^T\left(B_{\mathcal{N}_i}x_{\mathcal{N}_i} - y_i\right)\Big\}\\
\nonumber &=\displaystyle\sup_{\lambda,\mu}\sum_{i=1}^M\bigg[\inf_{x_i}\Big\{ 
f_i(x_i)+h_i(x_i)+x_i^T\big(
  A_{\mathcal{M}_i}^T\lambda_{\mathcal{M}_i}\\
\nonumber & \qquad\qquad~~+B_{\mathcal{M}_i}^T\mu_{\mathcal{M}_i}\big)\Big\}
-\lambda_i^Tb_i+\inf_{y_i}\Big\{
  g_i(y_i)-\mu_i^Ty_i\Big\}\bigg].
\end{align}
Introducing $F_i := f_i+h_i$ and $F := \sum_{i=1}^M F_i = \sum_{i=1}^M
f_i+h_i$, and noting the definition of
conjugate functions in the above
expression, we get that the dual problem can be written as
\begin{align}
\nonumber &\displaystyle\sup_{\lambda,\mu}\sum_{i=1}^M
\left\{-F_i^\star(-A_{\mathcal{M}_i}^T\lambda_{\mathcal{M}_i}-B_{\mathcal{M}_i}^T\mu_{\mathcal{M}_i})-\lambda_i^Tb_i-g_i^\star(\mu_i)\right\}\\
\label{eq:dualProb}&= \sup_{\lambda,\mu}\left\{-F^\star(-A^T\lambda-B^T\mu)-\lambda^Tb-g^\star(\mu)\right\}.
\end{align}
We further introduce $\nu =
(\lambda,\mu)\in\mathbb{R}^{m+p}$, $\nu_{\mathcal{M}_i} =
(\lambda_{\mathcal{M}_i},\mu_{\mathcal{M}_i})\in\mathbb{R}^{m_{\mathcal{M}_i}+p_{\mathcal{M}_i}}$,
$\nu_i =
(\lambda_i,\mu_i)\in\mathbb{R}^{m_i+p_i}$, $C =
[A^T~B^T]^T\in\mathbb{R}^{(m+p) \times n}$, $C_{\mathcal{M}_i}^T =
[A_{\mathcal{M}_i}^T~B_{\mathcal{M}_i}^T]\in\mathbb{R}^{n_i\times(m_{\mathcal{M}_i}+p_{\mathcal{M}_i})}$,
$c = (b,0)\in\mathbb{R}^{m+p}$, $c_i =
(b_i,0)\in\mathbb{R}^{m_i+p_i}$, and the following functions:
\begin{align}
d_i(\nu_{\mathcal{M}_i}) &:=
-F_{i}^\star(-C_{\mathcal{M}_i}^T\nu_{\mathcal{M}_i})-c_i^T\nu_i
\label{eq:localDualFcn}\\
\label{eq:dualFcn} d(\nu) &:= -F^\star(-C^T\nu)-c^T\nu
\end{align}
where $d = \sum_{i=1}^M d_i$. Using these definitions and notations we
arrive at the following dual problem:
\begin{align}
\nonumber &\displaystyle\sup_{\nu}\sum_{i=1}^M
\left\{d_i(\nu_{\mathcal{M}_i})-g_i^\star([0~I]\nu_i)\right\}\\
\label{eq:dualProbSimpl} &=\sup_{\nu}\left\{d(\nu)-g^\star([0~I]\nu)\right\}.
\end{align}
To evaluate $d_i$ or $d$ (or equivalently $F_i^\star$ or $F^\star$),
an optimization problem must be solved due to the definition of the
conjugate function. The minimands to these optimization problems are
defined by
\begin{align}
\label{eq:solInnerProbLocal}x_i^\star(\nu_{\mathcal{M}}) &:=
\arg\min_{x_i}\left\{f_i(x_i)+h_i(x_i)+\nu_{\mathcal{M}_i}^TC_{\mathcal{M}_i}x_i\right\},\\
\label{eq:solInnerProb}x^\star(\nu) &:= \arg\min_{x}\left\{f(x)+h(x)+\nu^TCx\right\}
\end{align}
since $F_i = f_i+h_i$ and $F=f+h$ respectively.
From Corollary~\ref{cor:strConvSmooth} we have that $d_i$ and $d$ are
differentiable with gradients
\begin{align*}
\nabla d_i(\nu_{\mathcal{M}_i}) &= C_{\mathcal{M}_i}x_i^\star(\nu_{\mathcal{M}_i})-\hat{c}_i,\\
\nabla d(\nu) &= Cx^\star(\nu)-c.
\end{align*}
respectively, where $\hat{c}_i=(0,\ldots,0,c_i,0,\ldots,0)$. Further, differentiation of the dual function w.r.t.
$\nu_i$ is given by
\begin{equation*}
\nabla_{\nu_i} d(\nu) = C_{\mathcal{N}_i}x_{\mathcal{N}_i}^\star(\nu_i)-c_i.
\end{equation*}
Corollary~\ref{cor:strConvSmooth} further implies that the gradients
to $d_i$ and $d$
are Lipschitz continuous with constants $L_i=
\|C_{\mathcal{M}_i}\|_2^2/\lambda_{\min}(H_i)$ and $L=
\|C\|_2^2/\lambda_{\min}(H)$ respectively. As previously discussed,
this is equivalent to the
existence of a quadratic lower bound given by
\eqref{eq:standardQuadBoundConcave} to the concave dual function, with
curvature $L_i$ and $L$ respectively.
In the following section we will show that the functions $d_i$ and $d$
defined in \eqref{eq:localDualFcn} and\eqref{eq:dualFcn}
respectively, satisfy the following tighter lower bounds
\begin{equation}
d(\nu_1)\geq d(\nu_2)+\langle
\nabla d(\nu_2),\nu_1-\nu_2\rangle-\tfrac{1}{2}\|\nu_1-\nu_2\|_{CH^{-1}C^T}^2
\label{eq:impDualBound}
\end{equation}
for all $\nu_1,\nu_2\in\mathbb{R}^{m+p}$ and
\begin{multline}
\label{eq:impLocalDualBound} d_i(\nu_{\mathcal{M}_i}^1)\geq d_i(\nu_{\mathcal{M}_i}^2)+\langle
\nabla
d_i(\nu_{\mathcal{M}_i}^2),\nu_{\mathcal{M}_i}^1-\nu_{\mathcal{M}_i}^2\rangle\\
-\tfrac{1}{2}\|\nu_{\mathcal{M}_i}^1-\nu_{\mathcal{M}_i}^2\|_{C_{\mathcal{M}_i}H_i^{-1}C_{\mathcal{M}_i}^T}
\end{multline}
for all
$\nu_{\mathcal{M}_i}^1,\nu_{\mathcal{M}_i}^2\in\mathbb{R}^{m_{\mathcal{M}_i}+p_{\mathcal{M}_i}}$
respectively.
We will also show that if the primal cost $f$ is a quadratic with positive
definite Hessian and $h$ is the indicator function for a closed, convex set,
no better quadratic lower bound exists.

\section{Dual function properties}

To show that the dual and local dual functions satisfy \eqref{eq:impDualBound} and
\eqref{eq:impLocalDualBound} respectively, some
preliminary results are needed. For notational convenience, we will
state the results for the function $d$ in the main parts of this
section. The corresponding  
results for the functions $d_i$ are given in the end.

In the following lemma we show that the distance in
$\|\cdot\|_H$-norm, where $H\in\mathbb{S}_{++}^n$ is the matrix
defining the strong convexity property of $f$, between any two points
$x^\star(\nu_1),x^\star(\nu_2)\in\mathbb{R}^n$ is upper bounded by
$\|\nu_1-\nu_2\|_{CH^{-1}C^T}$.
\begin{lem}
Suppose that Assumption~\ref{ass:probAss} holds and that $f$ is
strongly convex with matrix $H\in\mathbb{S}_{++}^n$. Then
\begin{equation*}
\|x^\star(\nu_1)-x^\star(\nu_2)\|_H\leq \|\nu_1-\nu_2\|_{CH^{-1}C^T}
\end{equation*}
for every $\nu_1,\nu_2\in\mathbb{R}^{(m+p)}$ where
$x^{\star}(\nu)$ is given by \eqref{eq:solInnerProb}, and
$C\in\mathbb{R}^{(m+p)\times n}$ is the equality
constraint matrix in
\eqref{eq:primProb}.
\label{lem:helpLemma}
\end{lem}
\begin{pf}
We first show that
\begin{multline}
\label{eq:lemmaHelpResult}\langle\nabla f(x^\star(\nu_1))-\nabla
f(x^\star(\nu_2)),x^\star(\nu_1)-x^\star(\nu_2)\rangle\leq
\\
\leq\langle C^T(\nu_1-\nu_2),x^\star(\nu_2)-x^\star(\nu_1)\rangle.
\end{multline}
First order optimality conditions for \eqref{eq:solInnerProb} using
$\nu_1$ and $\nu_2$ respectively are
\begin{align}
\label{eq:firstOrdOptIncl1}0&\in\nabla f(x^\star(\nu_1))+\partial
h(x^\star(\nu_1))+C^T\nu_1,\\
\label{eq:firstOrdOptIncl2}0&\in\nabla f(x^\star(\nu_2))+\partial
h(x^\star(\nu_2))+C^T\nu_2.
\end{align}
We denote by $\xi(x^\star(\nu_1))\in\partial h(x^\star(\nu_1))$
and $\xi(x^\star(\nu_2))\in\partial h(x^\star(\nu_2))$ the
sub-gradients that give equalities in \eqref{eq:firstOrdOptIncl1} and
\eqref{eq:firstOrdOptIncl2} respectively. This gives
\begin{align}
\label{eq:firstOrdOptEq1}0&=\nabla f(x^\star(\nu_1))+\xi(x^\star(\nu_1))+C^T\nu_1,\\
\label{eq:firstOrdOptEq2}0&=\nabla f(x^\star(\nu_2))+\xi(x^\star(\nu_2))+C^T\nu_2.
\end{align}
Taking the scalar product of \eqref{eq:firstOrdOptEq1} with
$x^\star(\nu_2)-x^\star(\nu_1)$ and the scalar product of
\eqref{eq:firstOrdOptEq2} with 
$x^\star(\nu_1)-x^\star(\nu_2)$, and summing the resulting expressions give
\begin{align*}
\langle\nabla f(x^\star(\nu_1))-\nabla
f(x^\star(\nu_2)),x^\star(\nu_1)-x^\star(\nu_2)\rangle+\qquad&\\
+\langle
C^T(\nu_1-\nu_2),x^\star(\nu_1)-x^\star(\nu_2)\rangle&=\\
=\langle\xi(x^\star(\nu_1))-\xi(x^\star(\nu_2)),x^\star(\nu_2)-x^\star(\nu_1)\rangle&\leq
0
\end{align*}
where the inequality holds since sub-differentials of proper,
closed, and convex functions are (maximal) monotone
mappings, see \cite[~\S 24]{Rockafellar}.
This implies that \eqref{eq:lemmaHelpResult} holds.

Further
\begin{align*}
  \|&x^\star(\nu_1)-x^\star(\nu_2)\|_H^2\leq\\
&\leq \langle\nabla f(x^\star(\nu_1))-\nabla
f(x^\star(\nu_2)),x^\star(\nu_1)-x^\star(\nu_2)\rangle\\
&\leq \langle
C^T(\nu_1-\nu_2),x^\star(\nu_2)-x^\star(\nu_1)\rangle
\\
&=\langle
H^{-1/2}C^T(\nu_1-\nu_2),H^{1/2}(x^\star(\nu_2)-x^\star(\nu_1))\rangle\\
&\leq \|H^{-1/2}C^T(\nu_1-\nu_2)\|_2\|x^\star(\nu_2)-x^\star(\nu_1)\|_H
\end{align*}
where the first inequality comes from
Lemma~\ref{lem:strConvCondition}, the second from \eqref{eq:lemmaHelpResult}, 
and the final inequality is due to Cauchy Schwarz.
This implies that
\begin{equation*}
\|x^\star(\nu_1)-x^\star(\nu_2)\|_H\leq \|\nu_1-\nu_2\|_{CH^{-1}C^T}
\end{equation*}
which concludes the proof.
\end{pf}

We are now ready to state the main theorem of this section.
\begin{thm}
Suppose that Assumption~\ref{ass:probAss} holds and that $f$ is
strongly convex with matrix $H\in\mathbb{S}_{++}^n$. The dual function $d$ defined in
\eqref{eq:dualFcn} is concave, differentiable and
satisfies 
\begin{equation}
d(\nu_1)\geq d(\nu_2)+\langle
\nabla d(\nu_2),\nu_1-\nu_2\rangle-\tfrac{1}{2}\|\nu_1-\nu_2\|_{\mathbf{L}}^2
\label{eq:dualQuadBound}
\end{equation}
for every $\nu_1,\nu_2\in\mathbb{R}^{m+p}$ and any
$\mathbf{L}\in\mathbb{S}_{+}^{m+p}$ such that $\mathbf{L}\succeq
CH^{-1}C^T$.
\label{thm:dualQuadBound}
\end{thm}
\begin{pf}
Concavity and differentiability is deduced from Danskin's Theorem, see
\cite[Proposition~B.25]{Bertsekas99}.

To show \eqref{eq:dualQuadBound}, we have for any
$\nu_1,\nu_2\in\mathbb{R}^{m}$ that
\begin{align*}
\langle\nabla d(\nu_1)-&\nabla
d(\nu_2),\nu_2-\nu_1\rangle=\\
&= \langle
Cx^\star(\nu_1)-c-Cx^\star(\nu_2)+c,\nu_2-\nu_1\rangle\\
&=\langle
x^\star(\nu_1)-x^\star(\nu_2),C^T(\nu_2-\nu_1)\rangle\\
&=\langle
x^\star(\nu_1)-x^\star(\nu_2),H^{-1}C^T(\nu_2-\nu_1)\rangle_H\\
&\leq \|
x^\star(\nu_1)-x^\star(\nu_2)\|_H\|H^{-1}C^T(\nu_2-\nu_1)\|_H\\
&\leq \|H^{-1}C^T(\nu_2-\nu_1)\|_H^2\\
&= (\nu_2-\nu_1)^TCH^{-1}C^T(\nu_2-\nu_1)\\
&= \|\nu_2-\nu_1\|_{CH^{-1}C^T}^2
\end{align*}
where the first inequality is due to Cauchy-Schwarz, and the
second is from Lemma~\ref{lem:helpLemma}.
Applying Corollary~\ref{cor:quadUBequivConcave} gives that
\eqref{eq:dualQuadBound} holds for every
$\nu_1,\nu_2\in\mathbb{R}^m$.
\end{pf}
\begin{cor}
The local dual functions $d_i$ defined in \eqref{eq:localDualFcn} are
concave, differentiable and satisfy
\begin{align*}
d_i(\nu_{\mathcal{M}_i}^1)&\geq d_i(\nu_{\mathcal{M}_i}^2)+\langle
\nabla
d_i(\nu_{\mathcal{M}_i}^2),\nu_{\mathcal{M}_i}^1-\nu_{\mathcal{M}_i}^2\rangle-\\
&\qquad\qquad\qquad\qquad\qquad-\tfrac{1}{2}\|\nu_{\mathcal{M}_i}^1-\nu_{\mathcal{M}_i}^2\|_{\mathbf{L}_{\mathcal{M}_i}}^2
\end{align*}
for all
$\nu_{\mathcal{M}_i}^1,\nu_{\mathcal{M}_i}^2\in\mathbb{R}^{m_{\mathcal{M}_i}+p_{\mathcal{M}_i}}$
and any
$\mathbf{L}_{\mathcal{M}_i}\in\mathbb{S}_{++}^{m_{\mathcal{M}_i}+p_{\mathcal{M}_i}}$ such that
$\mathbf{L}_{\mathcal{M}_i}\succeq C_{\mathcal{M}_i}H_i^{-1}C_{\mathcal{M}_i}^T$.
\label{cor:localDualQuadBound}
\end{cor}
\begin{pf}
The proof follows the same lines as the proof to
Theorem~\ref{thm:dualQuadBound}.
\end{pf}

Next, we show that if $f$ is a strongly convex quadratic function
and $h$ satisfies certain conditions, then
Theorem~\ref{thm:dualQuadBound} gives the best possible bound of the
form \eqref{eq:dualQuadBound}.
\begin{prop}
Assume that $f(x) = \tfrac{1}{2}x^THx+\zeta^Tx$ with
$H\in\mathbb{S}_{++}^n$ and $\zeta\in\mathbb{R}^n$ and that there exists a set
$\mathcal{X}\subseteq\mathbb{R}^n$ with non-empty interior on which
$h$ (besides being
proper, closed, and convex) is linear, i.e. $h(x) =
\xi_{\mathcal{X}}^Tx+\theta_\mathcal{X}$ for all $x\in\mathcal{X}$. Further,
assume that there exists $\widetilde{\nu}$ such that
$x^\star(\widetilde{\nu})\in{\rm{int}}(\mathcal{X})$.
Then for any matrix $\mathbf{L}\not\succeq
CH^{-1}C^T$, there exist $\nu_1$ and $\nu_2$ such that
\eqref{eq:dualQuadBound} does not hold.
\label{prp:dualQuadBoundTight}
\end{prop}
\begin{pf}
Since
$x^\star(\widetilde{\nu})\in{\rm{int}}(\mathcal{X})$ we get 
for all $\nu_{\epsilon}\in\mathcal{B}_{\epsilon}^{m+p}(0)$, where the radius
$\epsilon$ is small enough, that
$x^\star(\widetilde{\nu})-H^{-1}C^T\nu_{\epsilon}\in\mathcal{X}$.
Introducing $x_{\epsilon}=-H^{-1}C^T\nu_{\epsilon}$, we get from the
optimality conditions to \eqref{eq:solInnerProb}
(that specifies $x^{\star}(\nu)$) that
\begin{align*}
0&=Hx^\star(\widetilde{\nu})+\zeta+\xi_\mathcal{X}+C^T\widetilde{\nu}\\
&=H(x^\star(\widetilde{\nu})+x_{\epsilon})+\zeta+\xi_\mathcal{X}+C^T(\widetilde{\nu}+\nu_{\epsilon})\\
&=H(x^\star(\widetilde{\nu})+x_{\epsilon})+\zeta+h^{\prime}
(x^\star(\widetilde{\nu})+x_{\epsilon})+C^T(\widetilde{\nu}+\nu_{\epsilon})
\end{align*}
where $h^{\prime}(x^\star(\widetilde{\nu})\in\partial
h(x^\star(\widetilde{\nu})$ and
$x^\star(\widetilde{\nu})+x_{\epsilon}\in\mathcal{X}$ is used in
the last step. This implies that
$x^\star(\widetilde{\nu}+\nu_{\epsilon}) =
x^\star(\widetilde{\nu})+x_{\epsilon}$ and consequently that
$x^\star(\widetilde{\nu}+\nu_{\epsilon})\in\mathcal{X}$ for any
$\nu_{\epsilon}\in\mathcal{B}_{\epsilon}^{m+p}(0)$. Thus, for any
$\nu\in \widetilde{\nu}\oplus\mathcal{B}_{\epsilon}^{m+p}(0)$
we get
\begin{align*}
d(\nu) &= \min_{x}\tfrac{1}{2}x^THx+\zeta^Tx+h(x)+\nu^T(Cx-c)\\
&= \min_{x}\tfrac{1}{2}x^THx+\zeta^Tx+\xi_{\mathcal{X}}^Tx+\nu^T(Cx-c)\\
&=-\tfrac{1}{2}\nu^TCH^{-1}C^T\nu+\xi^T\nu+\theta
\end{align*}
where $\xi\in\mathbb{R}^n$ and $\theta\in\mathbb{R}$ collects the
linear and constant terms respectively. Since on the set
$\widetilde{\nu}\oplus\mathcal{B}_{\epsilon}^{m+p}(0)$,
$d$ is a quadratic with Hessian $CH^{-1}C^T$, it is straight-forward
to verify that \eqref{eq:dualQuadBound} holds
with equality for all
$\nu_1,\nu_2\in\widetilde{\nu}\oplus\mathcal{B}_{\epsilon}^{m+p}(0)$
if $\mathbf{L}=CH^{-1}C^T$. Thus, since
$\widetilde{\nu}\oplus\mathcal{B}_{\epsilon}^{m+p}(0)$ has non-empty
interior, we can for any matrix $\mathbf{L}\not\succeq
CH^{-1}C^T$ find
$\nu_1,\nu_2\in\widetilde{\nu}\oplus\mathcal{B}_{\epsilon}^{m+p}(0)$
such that 
\begin{equation*}
\|\nu_1-\nu_2\|_{CH^{-1}C^T}\geq \|\nu_1-\nu_2\|_{\mathbf{L}}.
\end{equation*}
This implies that for any $\mathbf{L}\not\succeq
CH^{-1}C^T$ there exist
$\nu_1,\nu_2\in\widetilde{\nu}\oplus\mathcal{B}_{\epsilon}^{m+p}(0)$
such that
\begin{align*}
d(\nu_1)&= d(\nu_2)+\langle \nabla
d(\nu_2),\nu_1-\nu_2\rangle-\tfrac{1}{2}\|\nu_1-\nu_2\|_{CH^{-1}C^T}\\
&\leq d(\nu_2)+\langle \nabla
d(\nu_2),\nu_1-\nu_2\rangle-\tfrac{1}{2}\|\nu_1-\nu_2\|_{\mathbf{L}}
\end{align*}
This concludes the proof.
\end{pf}
Proposition~\ref{prp:dualQuadBoundTight} shows that the bound in
Theorem~\ref{thm:dualQuadBound} is indeed the best obtainable bound of
the form \eqref{eq:dualQuadBound} if $f$ is a quadratic and $h$
specifies the stated assumptions. Examples of functions
that satisfy the assumptions on $h$ in
Proposition~\ref{prp:dualQuadBoundTight}  include linear functions, indicator
functions of closed convex constraint sets with non-empty interior, and the
1-norm.

The main results of this section,
Theorem~\ref{thm:dualQuadBound} and
Corollary~\ref{cor:localDualQuadBound}, provide a tighter
quadratic lower bound to the dual function 
compared to 
what has previously been presented in the literature, i.e. compared to
Proposition~\ref{prp:strConvSmooth} and
Corollary~\ref{cor:strConvSmooth}. These results are the key to constructing
more efficient distributed algorithms.

\section{Distributed optimization algorithm}

Dual decomposition methods often suffer from slow convergence 
properties, although the use of fast gradient methods have improved
the situation. In
this section, we describe one distributed and one parallel dual
decomposition method that improves the convergence of such
methods significantly. In the distributed algorithm, both
primal
and dual variables are updated distributively, while in the parallel
algorithm, the primal variables are updated in parallel and the dual
variables are updated centralized. We will show how
the results presented in Theorem~\ref{thm:dualQuadBound} and
Corollary~\ref{cor:localDualQuadBound} together with
generalized fast gradient methods, \cite{WangmengTIP2011}, are
combined to arrive at these algorithms and indicate why the improved
convergence is achieved.

Generalized fast gradient
methods can be applied to solve problems of the form
\begin{equation}
{\hbox{minimize }} \ell(x)+\psi(x)
\label{eq:minfP}
\end{equation}
where $x\in\mathbb{R}^n$,
$\psi~:~\mathbb{R}^n\to\mathbb{R}\cup\{\infty\}$ is proper, closed and
convex, $\ell~:~\mathbb{R}^n\to\mathbb{R}$ is convex, differentiable, and
satisfies 
\begin{equation}
\ell(x_1)\leq \ell(x_2)+\langle \nabla
\ell(x_2),x_1-x_2\rangle+\tfrac{1}{2}\|x_1-x_2\|_{\mathbf{L}}^2
\label{eq:quadUpperBound}
\end{equation}
for all $x_1,x_2\in\mathbb{R}^n$ and some
$\mathbf{L}\in\mathbb{S}_{++}^n$.
Before we state the algorithm, we define the generalized prox operator
\begin{equation}
{\rm{prox}}_{\psi}^{\mathbf{L}}(x) := \arg\min_y
\left\{\psi(y)+\tfrac{1}{2}\|y-x\|_{\mathbf{L}}^2\right\}
\label{eq:proxOpDef}
\end{equation}
and note that
\begin{align}
\label{eq:proxEquiv} &{\rm{prox}}_{\psi}^{\mathbf{L}}(x- \mathbf{L}^{-1}\nabla \ell(x)) \\
\nonumber &=\arg\min_y\left\{
\tfrac{1}{2}\|y-x+\mathbf{L}^{-1}\nabla \ell(x)\|_{\mathbf{L}}^2+\psi(y)\right\}\\
\nonumber &=\arg\min_y \left\{\ell(x)+\langle \nabla \ell(x),y-x\rangle+\tfrac{1}{2}\|y-x\|_{\mathbf{L}}^2+\psi(y)\right\}.
\end{align}
The generalized fast gradient method is stated below.

\medskip
\hrule
\smallskip
\begin{alg} ~\\\indent{\bf{Generalized fast gradient method}}
\smallskip\hrule\smallskip
\noindent Set: $y^1= x^0\in\mathbb{R}^n, t^1=1$\\
{\bf{For}} $k\geq 1$
\begin{itemize}
\item[] $x^{k} =
  {\rm{prox}}_{\psi}^{\mathbf{L}}(y^k-\mathbf{L}^{-1}\nabla \ell(y^k))$
\item[] $t^{k+1} = \frac{1+\sqrt{1+4(t^k)^2}}{2}$
\item[] $y^{k+1} = x^k+\left(\frac{t^k-1}{t^{k+1}}\right)(x^k-x^{k-1})$
\end{itemize}
\smallskip\hrule
\label{alg:GFGM}
\end{alg}
\medskip

The standard fast gradient method as presented in
\cite{BecTab_FISTA:2009} is obtained by setting $\mathbf{L} = LI$ in
Algorithm~\ref{alg:GFGM}, 
where $L$ is the Lipschitz constant to $\nabla \ell$.
The main step of the fast gradient method is to perform a
prox-step, i.e., to minimize \eqref{eq:proxEquiv} which can be seen as
an approximation of the function $\ell+\psi$. For the standard fast gradient
method, $\ell$ is approximated with a quadratic upper bound that has the
same curvature, described by $L$, in all directions. If this quadratic
upper bound is a bad approximation of the
function to be minimized, slow convergence is
expected. The generalization to allow for a matrix $\mathbf{L}$ in the
algorithm allows for quadratic upper bounds with different curvature
in different directions. This enables for quadratic upper bounds that
much better approximate the function $\ell$ and consequently gives
improved convergence properties.

The generalized fast gradient method has a convergence rate of (see
\cite{WangmengTIP2011}) 
\begin{equation}
  \ell_{\psi}(x^k)-\ell_{\psi}(x^\star)\leq \frac{2\|x^\star-x^0\|_{\mathbf{L}}^2}{(k+1)^2}
\label{eq:convRate}
\end{equation}
where $\ell_{\psi}:= \ell+\psi$. The convergence rate of the standard fast gradient
method as given in \cite{BecTab_FISTA:2009}, is obtained by setting
$\mathbf{L}=LI$ in \eqref{eq:convRate}.

The objective here is to apply the generalized fast gradient method to
solve the dual problem \eqref{eq:dualProb}. By introducing
$\widetilde{g}(\nu)=g^\star([0~I]\nu)$,
the dual problem \eqref{eq:dualProb} can be expressed $\max_{\nu}
d(\nu)-\widetilde{g}(\nu)$, where
$d$ is defined in \eqref{eq:dualFcn}.
As shown in Theorem~\ref{thm:dualQuadBound}, the function $-d$ satisfies the
properties required to apply generalized fast gradient methods. Namely
that \eqref{eq:quadUpperBound} holds for
any $\mathbf{L}\in\mathbb{S}_{+}^{m+p}$ such that $\mathbf{L}\succeq CH^{-1}C^T$.
Further, since $g$ is a
closed, proper, and convex function so is $g^\star$, see \cite[Theorem
12.2]{Rockafellar}, and by \cite[Theorem
5.7]{Rockafellar} so is $\widetilde{g}$. This implies
that generalized fast gradient methods, i.e. Algorithm~\ref{alg:GFGM},
can be used to solve the dual problem \eqref{eq:dualProb}.
We set $-d = \ell$ and $\widetilde{g}=\psi$, and restrict
$\mathbf{L}={\rm{blkdiag}}(\mathbf{L}_{\lambda},\mathbf{L}_{\mu})$ to get the following
algorithm.

\medskip
\hrule
\smallskip
\begin{alg}~\\\noindent {\bf{Generalized fast dual gradient method}}
\smallskip\hrule\smallskip
\noindent Set: $z^1= \lambda^0\in\mathbb{R}^m, v^1=\mu^0\in\mathbb{R}^p, t^1=1$\\
{\bf{For}} $k\geq 1$
\begin{itemize}
\item[] $y^k = \arg\min_{x}\left\{f(x)+h(x)+(z^k)^TAx+(v^k)^TBx\right\}$
\item[] $\lambda^{k} = z^k+\mathbf{L}_{\lambda}^{-1}(Ay^k-b)$
\item[] $\mu^{k} =
  {\rm{prox}}_{g^\star}^{\mathbf{L}_{\mu}}(v^k+\mathbf{L}_{\mu}^{-1}By^k)$
\item[] $t^{k+1} = \frac{1+\sqrt{1+4(t^k)^2}}{2}$
\item[] $z^{k+1} = \lambda^k+\left(\frac{t^k-1}{t^{k+1}}\right)(\lambda^k-\lambda^{k-1})$
\item[] $v^{k+1} =
  \mu^k+\left(\frac{t^k-1}{t^{k+1}}\right)(\mu^k-\mu^{k-1})$
\end{itemize}
\smallskip\hrule
\label{alg:GFDGM}
\end{alg}
\medskip
where $y^k$ is the primal variable at iteration $k$ that is used to
help compute the gradient $\nabla
d(\nu^k)$ where $\nu^k =
(z^k,v^k)$. To arrive
at the $\lambda^k$ and $\mu^k$ iterations, we let $\xi^k =
(\lambda^k,\mu^k)$,
and note that
\begin{align}
\label{eq:proxDecomposition} \xi^k&={\rm{prox}}_{\widetilde{g}}^{\mathbf{L}}(\nu^k+\mathbf{L}^{-1}\nabla d(\nu^k))\\
\nonumber &=\arg\min_{\nu} \left\{\tfrac{1}{2}\|\nu-\nu^k-\mathbf{L}^{-1}\nabla
d(\nu^k)\|_{\mathbf{L}}^2+g^\star([0
~I]\nu)\right\}\\
\nonumber &= \left[\begin{array}{l}
\arg\min_{z} \left\{\tfrac{1}{2}\|z-z^k-\mathbf{L}_{\lambda}^{-1}\nabla_{z}
d(\nu^k)\|_{\mathbf{L}_{\lambda}}^2\right\}\\
\arg\min_{v} \big\{\tfrac{1}{2}\|v-v^k-\mathbf{L}_{\mu}^{-1}\nabla_{v}
  d(\nu^k)\|_{\mathbf{L}_{\mu}}^2+g^\star(v)\big\}
\end{array}\right]\\
\nonumber &=\left[\begin{array}{l}
z^k+\mathbf{L}_{\lambda}^{-1}(Ay^k-b)\\
{\rm{prox}}_{g^\star}^{\mathbf{L}_{\mu}}(v^k+\mathbf{L}_{\mu}^{-1}By^k)
\end{array}\right].
\end{align}

When solving separable
problems of the form \eqref{eq:primProbSep}, Algorithm~\ref{alg:GFDGM} can be implemented
in distributed fashion by restricting
$\mathbf{L}_{\lambda}\in\mathbb{S}_{++}^{m}$ and $\mathbf{L}_{\mu}\in\mathbb{S}_{++}^{p}$
to be block diagonal, i.e. 
$\mathbf{L}_{\lambda}={\rm{blkdiag}}(\mathbf{L}_{\lambda
  1},\ldots,\mathbf{L}_{\lambda M})$ and 
$\mathbf{L}_{\mu}={\rm{blkdiag}}(\mathbf{L}_{\mu
  1},\ldots,\mathbf{L}_{\mu M})$ and 
where $\mathbf{L}_{\lambda i}\in\mathbb{S}_{++}^{m_i}$ and
$\mathbf{L}_{\mu i}\in\mathbb{S}_{++}^{p_i}$. The distributed implementation
is presented next.

\medskip
\hrule
\smallskip
\begin{alg}\label{alg_dist_acc_prox} ~\\
\textbf{Distributed generalized fast dual gradient method}
\vspace{-2mm}\hrule\smallskip
\noindent Initialize
$z_i^1=\lambda_i^{0}\in\mathbb{R}^{m_i}, v_i^1=\mu_i^0\in\mathbb{R}^{p_i}, t^1 = 1$.\\
In every node, $i=\{1,\ldots,M\}$, do the following steps\\
{\bf{For}} $k \geq 1$
\begin{enumerate}
\item Send $z_i^{k},v_i^k$ to each $j \in
   \mathcal{N}_i$,
\item[] receive $z_j^{k},v_j^k$ from each $j\in 
   \mathcal{M}_i$
\item Form $z_{\mathcal{M}_i} = (\ldots,z_j^k,\ldots)$ with all $j\in\mathcal{M}_i$
\item Form $v_{\mathcal{M}_i} = (\ldots,v_j^k,\ldots)$ with all
  $j\in\mathcal{M}_i$
 \item Update local primal variables according to
\item[] $\displaystyle y_i^k =\arg\min_x
\left\{f_i(x)+h_i(x)+x_i^T\left(A_{\mathcal{M}_i}^Tz_{\mathcal{M}_i}^k+B_{\mathcal{M}_i}^Tv_{\mathcal{M}_i}^k\right)\right\}$
 \item Send $y_i^{k}$ to each $j \in \mathcal{M}_i$, receive
   $y_j^{k}$ from each $j \in \mathcal{N}_i$
 \item Form $y_{\mathcal{N}_i} = (\ldots,y_j^k,\ldots)$ with all $j\in\mathcal{N}_i$
 \item Update local dual variables according to
   \begin{itemize}
\item[] $\lambda_i^k =
  z_i^k+\mathbf{L}_{\lambda i}^{-1}(A_{\mathcal{N}_i}y_{\mathcal{N}_i}^k-b_i)$
\item[] $\mu_i^k = {\rm{prox}}_{g_i^\star}^{\mathbf{L}_{\mu
      i}}(v_i^k+\mathbf{L}_{\mu i}^{-1}B_{\mathcal{N}_i}y_{\mathcal{N}_i}^k)$
\item[] $t^{k+1} = \frac{1+\sqrt{1+4(t^k)^2}}{2}$
\item[] $z_i^{k+1} = \lambda_i^k+\left(\frac{t^k-1}{t^{k+1}}\right)(\lambda_i^k-\lambda_i^{k-1})$
\item[] $v_i^{k+1} =
  \mu_i^k+\left(\frac{t^k-1}{t^{k+1}}\right)(\mu_i^k-\mu_i^{k-1})$
\end{itemize}
\end{enumerate}
\smallskip\hrule
\label{alg:DGFDGM}
\end{alg}
\medskip

In this distributed algorithm, both the primal and dual variables are
updated in distributed fashion. When solving optimization problems
\eqref{eq:primProbSep} with all $g_i=0$, 
Algorithm~\ref{alg:DGFDGM} can be efficiently implemented in parallel
fashion in which the primal variables are updated in
parallel, while the dual variables are updated in a central unit.
A parallel implementation relaxes the block-diagonal requirement on
$\mathbf{L}$ which can 
give a considerably improved convergence rate.

\medskip
\hrule
\smallskip
\begin{alg}\label{alg_dist_acc_prox} ~\\
\textbf{Parallel generalized fast dual gradient method}
\smallskip\hrule\smallskip
\noindent Initialize
$z^1=(z_1^1,\ldots,z_M^k)=\lambda^{0}\in\mathbb{R}^{m}, t^1 = 1$.\\
{\bf{For}} $k \geq 1$
\begin{enumerate}
\item Form $z_{\mathcal{M}_i}^k = (\ldots,z_j^k,\ldots)$ with all $j\in\mathcal{M}_i$
\item Send $z_{\mathcal{M}_i}^{k}$ to each node $j \in \{1,\ldots,M\}$
\item Update local primal variables according to
\begin{itemize}
\item[] $\displaystyle y_i^k =\arg\min_x\left\{
f_i(x)+h_i(x)+x_i^TA_{\mathcal{M}_i}^Tz_{\mathcal{M}_i}^k\right\}$
\end{itemize}
 \item Receive $y_i^{k}$ from each node $j \in \{1,\ldots,M\}$
\item Form $y^k = (y_1^k,\ldots,y_M^k)$
 \item Update dual variables according to
   \begin{itemize}
\item[] $\lambda^k =
  z^k+\mathbf{L}^{-1}(Ay^k-b)$
\item[] $t^{k+1} = \frac{1+\sqrt{1+4(t^k)^2}}{2}$
\item[] $z^{k+1} = \lambda^k+\left(\frac{t^k-1}{t^{k+1}}\right)(\lambda^k-\lambda^{k-1})$
\end{itemize}
\end{enumerate}
\smallskip\hrule
\label{alg:PGFDGM}
\end{alg}
\medskip

The matrix $\mathbf{L}\in\mathbb{S}_{++}^m$ in
Algorithm~\ref{alg:PGFDGM} must satisfy
$\mathbf{L}\succeq AH^{-1}A^T$ (since $p=0$ and $B=0$ due to the assumption that
$g=\sum_ig_i=0$). Since $A$ by assumption is sparse and
has full row rank and $H$ is block-diagonal, we can choose $\mathbf{L}
= AH^{-1}A^T$. This gives
the tightest possible quadratic upper bound to the function $-d$, i.e.
we get a good approximation of $-d$ in the algorithm. When
implementing the algorithm, the inverse $\mathbf{L}^{-1}$ is obviously not computed in each
iteration. Rather, a sparse Cholesky or LDL-factorization of the
matrix $AH^{-1}A^T$ is computed offline and the factors are
stored for online use. Such sparse Cholesky and LDL-factorizations can
be computed for very large matrices.
This implies that inversion of the
$\mathbf{L}$-matrix in the algorithm reduces to one forward and one
backward solve for the sparse triangular factor and its transpose.
This can be very efficiently implemented.

In the following proposition we state the convergence rate properties
of Algorithm~\ref{alg:DGFDGM} and Algorithm~\ref{alg:PGFDGM}.
\begin{prop}
Suppose that Assumption~\ref{ass:probAss} holds. If, independent of
structure, 
$\mathbf{L}\succeq CH^{-1}C^T$ and $\mathbf{L}\succeq AH^{-1}A^T$ in
Algorithm~\ref{alg:DGFDGM} and 
Algorithm~\ref{alg:PGFDGM} respectively. Then Algorithm~\ref{alg:DGFDGM} and
Algorithm~\ref{alg:PGFDGM}
converges with the rate
\begin{equation}
D(\nu^\star)-D(\nu^k)
\leq \frac{2 \left\|
\nu^\star - \nu^0\right\|_{\mathbf{L}}^2}{(k+1)^2}, \forall k\geq 1
\label{eq:dualConvRate}
\end{equation}
where $D=d-\widetilde{g}$ and $k$ is the iteration number.
\label{prp:convRate}
\end{prop}
\begin{pf}
Algorithm~\ref{alg:DGFDGM} is a distributed and
Algorithm~\ref{alg:PGFDGM} is a parallel implementation of
Algorithm~\ref{alg:GFDGM}. They therefore share 
the same convergence rate properties. Algorithm~\ref{alg:GFDGM} is
Algorithm~\ref{alg:GFGM} applied to solve
the dual problem \eqref{eq:dualProb}.  
The convergence rate of Algorithm~\ref{alg:GFGM} is
given by
\eqref{eq:convRate} provided that the function to be minimized a sum
of one convex, differentiable function that satisfies \eqref{eq:quadUpperBound} and
one closed, proper, and convex function, see
\cite{WangmengTIP2011}. The discussion preceding the presentation of
Algorithm~\ref{alg:GFDGM} shows that the dual function to be optimized satisfies these
properties for any $\mathbf{L}\succeq CH^{-1}C^T$. This proves the
convergence rate for Algorithm~\ref{alg:DGFDGM}. Further for
Algorithm~\ref{alg:PGFDGM}, $g=0$, which implies $B=0$ and $C=A$. This
gives the conditions for Algorithm~\ref{alg:PGFDGM} and concludes the proof.
\end{pf}

\begin{rem}
By forming a specific running average of previous primal variables, it
is possible to prove a $O(1/k)$ convergence rate for the distance to the
primal variable optimum and a $O(1/k^2)$ convergence rate for
the worst case primal infeasibility, see \cite{PatrinosTAC2013}.
\end{rem}

For some choices of conjugate functions $g^\star$ and $g_i^{\star}$,
${\rm{prox}}_{g^\star}^{\mathbf{L}_{\mu}}(x)$ and 
${\rm{prox}}_{g_i^\star}^{\mathbf{L}_{\mu}}(x_i)$ in
Algorithm~\ref{alg:DGFDGM} can be difficult 
to evaluate. For
standard prox operators (given by ${\rm{prox}}_{g^\star}^{I}(x))$,
Moreau decomposition \cite[Theorem 31.5]{Rockafellar} states that
\begin{equation*}
{\rm{prox}}_{g^\star}^{I}(x)+{\rm{prox}}_{g}^{I}(x)=x.
\end{equation*}
In the following proposition, we will generalize this result to hold for
the generalized prox-operator used here.
\begin{prop}
Assume that $g~:~\mathbb{R}^n\to\mathbb{R}$ is a proper, closed, and
convex function. Then
\begin{equation*}
{\rm{prox}}_{g^\star}^{\mathbf{L}}(x)+\mathbf{L}^{-1}{\rm{prox}}_{g}^{\mathbf{L}^{-1}}(\mathbf{L}x)=x
\end{equation*}
for every $x\in\mathbb{R}^n$ and any $\mathbf{L}\in\mathbb{S}_{++}^n$.
\label{prp:genMoreau}
\end{prop}
\begin{pf}
Optimality conditions for the prox operator \eqref{eq:proxOpDef} give that
$y={\rm{prox}}_{g^\star}^{\mathbf{L}}(x)$ if and only if
\begin{equation*}
0\in \partial g^\star(y)+\mathbf{L}(y-x)
\end{equation*}
Introducing $v=\mathbf{L}(x-y)$ gives $v\in\partial g^\star(y)$ which
is equivalent to $y\in\partial g(v)$ \cite[Corollary 23.5.1]{Rockafellar}. Since $y =
x-\mathbf{L}^{-1}v$ we have
\begin{equation*}
0\in \partial g(v)+(\mathbf{L}^{-1}v-x)
\end{equation*}
which is the optimality condition for $v =
{\rm{prox}}_{g}^{\mathbf{L}^{-1}}(\mathbf{L}x)$. This concludes the proof.
\end{pf}
\begin{rem}
If $g=I_{\mathcal{X}}$ where $I_{\mathcal{X}}$ is the indicator
function, then $g^\star$ is the support function. Evaluating the prox
operator \eqref{eq:proxOpDef}
with $g^\star$ being a support function is difficult. However, through
Proposition~\ref{prp:genMoreau}, this can be rewritten to only require
the a projection operation onto the set $\mathcal{X}$.
If $\mathcal{X}$ is a box constraint and $\mathbf{L}$ is diagonal, then
the projection becomes a max-operation and hence very cheap to implement.
\end{rem}
\begin{rem}
Due to error accumulation of the fast gradient method, see \cite{Devolder}, the inner
minimizations, i.e. the $y_i^k$-updates, should be solved to high accuracy.
\end{rem}

We have shown how the $\mathbf{L}$-matrix should be chosen in the
parallel Algorithm~\ref{alg:PGFDGM}. However, we have not discussed
how to choose the block-diagonal $\mathbf{L}$-matrix used in  
Algorithm~\ref{alg:DGFDGM}.
This is the topic of the following section.

\section{Choosing the $\mathbf{L}$-matrix}

The (optimal) step-size selection in standard fast dual gradient methods relies on
computing a (tight) Lipschitz constant to the dual gradient. This
Lipschitz constant is usually computed by taking the Euclidean operator norm of
the equality constraint matrix $A$
(see Corollary~\ref{cor:strConvSmooth}). This requires centralized
computations. In this section we will extend a recent result in
\cite{BeckEtal2013} to allow for distributed selection of the
$\mathbf{L}$-matrix that is used in Algorithm~\ref{alg:DGFDGM}.

The $\mathbf{L}$-matrix in Algorithm~\ref{alg:DGFDGM} should be block
diagonal, i.e.
$\mathbf{L}={\rm{blkdiag}}(\mathbf{L}_1,\ldots,\mathbf{L}_M)$ to
facilitate a distributed implementation, and that it should satisfy
$\mathbf{L} \succeq CH^{-1}C^T$ to guarantee convergence of the
algorithm. We will see that
Corollary~\ref{cor:localDualQuadBound} can be used to compute a matrix
$\mathbf{L}$ that satisfies these requirements, using local computations
and neighboring communication only.
From Corollary~\ref{cor:localDualQuadBound} we have that any matrix
$\mathbf{L}_{\mathcal{M}_i}\in\mathbb{S}_{++}^{m_{\mathcal{M}_i}+p_{\mathcal{M}_i}}$
that describe a quadratic upper bound to the local dual functions
$d_i$ must satisfy $\mathbf{L}_{\mathcal{M}_i}\succeq
C_{\mathcal{M}_i}H_i^{-1}C_{\mathcal{M}_i}^T$. To allow for a
distributed implementation, we further
restrict $\mathbf{L}_{\mathcal{M}_i}$ to be block-diagonal, i.e. if
$\mathcal{M}_1 = \{1,4,6\}$
then $\mathbf{L}_{\mathcal{M}_1} =
{\rm{blkdiag}}(\mathbf{L}_{\mathcal{M}_1,1},\mathbf{L}_{\mathcal{M}_1,4},\mathbf{L}_{\mathcal{M}_1,6})$ where
$\mathbf{L}_{\mathcal{M}_i,j}\in\mathbb{S}_{++}^{m_j+p_j}$. These
restrictions on the local matrices $\mathbf{L}_{\mathcal{M}_i}$ are
summarized in the following set notation
\begin{align*}
\mathcal{L}_{\mathcal{M}_i} =
\big\{\mathbf{L}_{\mathcal{M}_i}\in&\mathbb{S}_{++}^{m_{\mathcal{M}_i}+p_{\mathcal{M}_i}}~|~\mathbf{L}_{\mathcal{M}_i}\succeq
C_{\mathcal{M}_i}H_i^{-1}C_{\mathcal{M}_i}^T,\\
& \qquad\quad\mathbf{L}_{\mathcal{M}_i} =
{\rm{blkdiag}}(\ldots,\mathbf{L}_{\mathcal{M}_i,j},\ldots)\\
& \qquad\quad{\hbox{with all }} j\in\mathcal{M}_i,
\mathbf{L}_{\mathcal{M}_i,j}\in\mathbb{S}_{++}^{m_j+p_j}\big\}.
\end{align*}
Using this set notation, we propose the following distributed
initialization procedure for Algorithm~\ref{alg:DGFDGM}.

\medskip
\hrule
\smallskip
\begin{alg}~\\\indent {\bf{Distributed initialization of Algorithm~\ref{alg:DGFDGM}}}
\smallskip\hrule\smallskip
\noindent For each $i\in \{1,\ldots,M\}$ \\
{\bf{Do}} 
\begin{enumerate}
\item Choose $\mathbf{L}_{\mathcal{M}_i}={\rm{blkdiag}}(\ldots,\mathbf{L}_{\mathcal{M}_i,j},\ldots)\in\mathcal{L}_{\mathcal{M}_i}$
\item Send $\mathbf{L}_{\mathcal{M}_i,j}$ to all
  $j\in\mathcal{M}_i$\\
Receive $\mathbf{L}_{\mathcal{M}_j,i}$ from all
  $j\in\mathcal{N}_i$
\item Compute $\mathbf{L}_i = \sum_{j\in\mathcal{N}_i}
  \mathbf{L}_{\mathcal{M}_j,i}$
\end{enumerate}
\smallskip\hrule
\label{alg:distrInit}
\end{alg}
\medskip

From this initialization we get local $\mathbf{L}_i$-matrices that
are used in each local node $i$ and in all iterations of
Algorithm~\ref{alg:DGFDGM}. In the following proposition we show that
Algorithm~\ref{alg:DGFDGM} converges
with the rate \eqref{eq:convRate} when initialized using Algorithm~\ref{alg:distrInit}.
\begin{prop}
Suppose that Assumption~\ref{ass:probAss} holds. If
$\mathbf{L}_i\in\mathbb{S}_{++}^{m_i+p_i}$ is computed using
Algorithm~\ref{alg:distrInit}. Then Algorithm~\ref{alg:DGFDGM}
converges with the rate \eqref{eq:convRate} when solving problems of the form
\eqref{eq:primProbSep}.
\label{prp:convRate2}
\end{prop}
\begin{pf}
For any $\nu=[\nu_1^T,\ldots,\nu_M^T]^T\in\mathbb{R}^{m+p}$,
and due to the notation
$\nu_{\mathcal{M}_i}\in\mathbb{R}^{m_{\mathcal{M}_i}+p_{\mathcal{M}_i}}$,
we get
\begin{align*}
\|\nu\|_{\mathbf{L}}^2&=\sum_{i=1}^M\|\nu_i\|_{\mathbf{L}_i}^2=\sum_{i=1}^M\sum_{j\in\mathcal{N}_i}\|\nu_i\|_{\mathbf{L}_{\mathcal{M}_j,i}}^2=\\
&=\sum_{i=1}^M\sum_{j\in\mathcal{M}_i}\|\nu_j\|_{\mathbf{L}_{\mathcal{M}_i,j}}^2=\sum_{i=1}^M
\|\nu_{\mathcal{M}_i}\|_{\mathbf{L}_{\mathcal{M}_i}}^2\geq\\
&\geq \sum_{i=1}^M
\|\nu_{\mathcal{M}_i}\|_{C_{\mathcal{M}_i}H_i^{-1}C_{\mathcal{M}_i}^T}^2
=\|\nu\|_{CH^{-1}C^T}^2
\end{align*}
which is equivalent to $\mathbf{L}\succeq CH^{-1}C^T$. Applying Proposition~\ref{prp:convRate}
completes the proof.
\end{pf}

The first step in the distributed initialization algorithm is still
not completely specified, i.e., we have not yet discussed how to choose
$\mathbf{L}_{\mathcal{M}_i}$. Since the primary application for our method
is distributed model predictive control (DMPC) in which similar optimization
problems are solved repeatedly online, much offline computational
effort can be devoted to ease the online computational burden. In the DMPC context, we
propose to solve the following local
optimization problem in step 1 and for each $i\in\{1,\ldots,M\}$:
\begin{equation}
\label{eq:findLocalL}
\begin{tabular}{ll}
minimize
& ${\hbox{tr }} \mathbf{L}_{\mathcal{M}_i}$\\
subject to & $\mathbf{L}_{\mathcal{M}_i} =
{\rm{blkdiag}}(\ldots,\mathbf{L}_{\mathcal{M}_i,j},\ldots)\in\mathcal{L}_{\mathcal{M}_i}$.
\end{tabular}
\end{equation}
This is a convex semi-definite program (SDP) that can readily be solved
using standard software. Another option in choosing
$\mathbf{L}_{\mathcal{M}_i}$ is to minimize the condition number of
$C_{\mathcal{M}_i}H_i^{-1}C_{\mathcal{M}_i}^T$, subject to structural
constraints. However, the condition number
is defined only if $C_{\mathcal{M}_i}$ has full row rank. For the case
of $C_{\mathcal{M}_i}$ having full column rank, the ratio between the
largest and smallest non-zero eigenvalues can be minimized. This is
achieved by minimizing the condition number of
$H_i^{-1/2}C_{\mathcal{M}_i}^TC_{\mathcal{M}_i}H_i^{-1/2}$. See
\cite[Section 3.1]{BoydLMI} and \cite[Section 6]{gis2014AutPart2} for
more on minimization of condition numbers and the ratio between the
largest and smallest eigenvalues of a symmetric positive semi-definite
matrix.

\section{Distributed model predictive control}

Distributed model predictive control (DMPC) is a distributed
optimization-based control scheme applied to control systems
consisting of several subsystems that have a sparse dynamic interaction
structure. The local dynamics are described by
\begin{align*}
x_{i}(t+1) &= \sum_{j\in\mathcal{N}_i}\Phi_{ij} x_{j}(t)+\Gamma_{ij}
u_{j}(t),& x_{i}(0)&=\bar{x}_i
\end{align*}
for all $i\in\{1,\ldots,M\}$, where $x_i\in\mathbb{R}^{n_{x_i}}$,
$u_i\in\mathbb{R}^{n_{u_i}}$, $\Phi_{ij}\in\mathbb{R}^{n_{x_i}\times
  n_{x_j}}$, $\Gamma_{ij}\in\mathbb{R}^{n_{x_i}\times
  n_{u_j}}$, and $\bar{x}_i\in\mathbb{R}^{n_{x_i}}$ is a measurement
of the current state. In DMPC, it is common 
to have local state and control
constraint sets $x_i\in\mathcal{X}_i$, $u_i\in\mathcal{U}_i$,
where $\mathcal{X}_i$ and $\mathcal{U}_i$ are non-empty, closed, and
convex sets. The cost function is usually chosen as the following sum over a
horizon $N$ 
\begin{equation*}
\sum_{i=1}^M\left(\sum_{t=0}^{N-1}\frac{1}{2}\begin{bmatrix}
x_i(t)\\
u_i(t)
\end{bmatrix}^T\begin{bmatrix}
Q_i & 0\\
0 & R_i
\end{bmatrix}
\begin{bmatrix}
x_i(t)\\
u_i(t)
\end{bmatrix}
\right)+\frac{1}{2}\|x_i(N)\|_{Q_{i,f}}^2
\end{equation*}
where $Q_i\in\mathbb{S}_{++}^{n_{x_i}}$,
$R_i\in\mathbb{S}_{++}^{n_{u_i}}$, 
and $Q_{i,f}\in\mathbb{S}_{++}^{n_{x_i}}$. By stacking the local state and control vectors into
$y_i=[x_i(0)^T,\ldots,x_i(N)^T,u_i(0)^T,\ldots,u_i(N-1)^T]^T$ 
we get an optimization problem of the form
\begin{equation}
\label{eq:DMPCoptProb}
\begin{tabular}{ll}
minimize & $\displaystyle\sum_{i=1}^M f_i(y_i)+h_i(y_i)$\\
subject to & $\displaystyle\sum_{j\in\mathcal{N}_i} A_{ij}y_j=b_i\bar{x}_i$
\end{tabular}
\end{equation}
where $f_i(y_i) = \frac{1}{2}y_i^TH_iy_i$, $h_i(y_i) =
I_{\mathcal{Y}_i}(y_i)$, and $H_i$, $\mathcal{Y}_i$, $A_{ij}$, and
$b_i$ are structured according to the stacked vector $y_i$. The
optimization problem \eqref{eq:DMPCoptProb} is structured as
\eqref{eq:primProbSep} and can therefore be
solved in distributed fashion using Algorithm~\ref{alg:DGFDGM} or in
parallel fashion using Algorithm~\ref{alg:PGFDGM}.

DMPC-problem with coupled linear inequality constraints
also fit into the framework presented in this paper. The corresponding
optimization problem becomes 
\begin{equation}
\label{eq:DMPCoptProbGen} 
\begin{tabular}{ll}
minimize & $\displaystyle\sum_{i=1}^M f_i(y_i)+h_i(y_i)+g_i(z_i)$\\
subject to & $\displaystyle\sum_{j\in\mathcal{N}_i} A_{ij}y_j=b_i\bar{x}_i$\\
 & $\displaystyle\sum_{j\in\mathcal{N}_i} B_{ij}y_j=z_i$
\end{tabular}
\end{equation}
where again $f_i(y_i) = \frac{1}{2}y_i^TH_iy_i$, $h_i(y_i) =
I_{\mathcal{Y}_i}(y_i)$, and $H_i$, $\mathcal{Y}_i$, $A_{ij}$, and
$b_i$ are structured according to the stacked vector $y_i$. The functions $g_i$
are the indicator functions for the coupled inequality constraints,
and the additional equality constraints
$\sum_{j\in\mathcal{N}_i} B_{ij}y_j=z_i$ describes the coupling.

Some formulations in the literature also use a coupled 1-norm cost for
reference tracking purposes, as in \cite{DoanHPV2013}. This also naturally fits
into the developed framework by letting $g_i(z_i)=\|z_i\|_1$ in
\eqref{eq:DMPCoptProbGen}.

We conclude this section with a remark on reconfigurability of the
proposed scheme in the DMPC context.
\begin{rem}
Due to the distributed structure of
the initialization procedure in Algorithm~\ref{alg:distrInit}, the DMPC scheme enjoys distributed
reconfiguration, commonly referred to as \emph{plug-and-play}.
Distributed reconfiguration or plug-and-play refers to the feature
that if an additional subsystem is connected to (or removed from) the system, the only
updates needed in the algorithm involve computations in the direct
neighborhood of the added (removed) subsystem. This is the case for
Algorithm~\ref{alg:distrInit} since if a reconfiguration is needed due
to addition or removal of subsystem $i$,
only subsystems $j\in\mathcal{M}_i$ need to be invoked
for the reconfiguration.
\end{rem}

\section{Numerical example}

\begin{table*}
\centering
\caption{Numerical evaluation between Algorithm~\ref{alg:DGFDGM},
  Algorithm~\ref{alg:PGFDGM}, fast dual decomposition, and the dual
  Newton CG method in \cite{dualNewtonCG}.}
\begin{tabular}{llcrrrrrrrr}
& &  &  \multicolumn{4}{c}{$\#$
  communication rounds} & \\
& &  &  \multicolumn{2}{c}{
  local} &\multicolumn{2}{c}{
  global} & avg. exec. time \\
Algorithm & Parameters & $\#$ ss/vars./constr. & avg. & max & avg. &
max &  12 cores [mm:ss.s] \\
\hline
Algorithm~\ref{alg:PGFDGM} & $\mathbf{L} = AH^{-1}A^T$  & 500/87k/246k 
&-&-& 16.2 & 118 & 2.3\\
Algorithm~\ref{alg:DGFDGM} & $\mathbf{L}$ computed using
Alg.~\ref{alg:distrInit}   & 500/87k/246k 
& 523.7 & 774 &-&-& 3.2 \\
Algorithm~\ref{alg:DGFDGM} & $\mathbf{L}= \|AH^{-1}A^T\|_2I$    & 500/87k/246k
 & 6114.7 & 6556 &-&-& 32.4 \\
Algorithm~\ref{alg:DGFDGM} & $\mathbf{L}= \|AH^{-1}A^T\|_1I$    & 500/87k/246k
 & 9923.2 & 10622 &- &-& 52.7 \\
\cite{dualNewtonCG} & $\epsilon_i=10^{-4}, \mu=0.8,\sigma=0.3$    & 500/87k/246k
 & 6661.1 & 28868 & 4082.6 & 17694 & 2:06.0 \\
\hline
Algorithm~\ref{alg:PGFDGM} & $\mathbf{L} = AH^{-1}A^T$  & 2000/351k/993k
 & -&-& 4.5 & 12 & 7.9 \\
Algorithm~\ref{alg:DGFDGM} & $\mathbf{L}$ blk-diag comp. fr.  & 2000/351k/993k
 & 356.8 & 652 &-&-& 15.6 \\
Algorithm~\ref{alg:DGFDGM} & $\mathbf{L}= \|AH^{-1}A^T\|_2I$   & 2000/351k/993k
 & 4474.9 & 4608 &-&-&  2:09:9 \\
Algorithm~\ref{alg:DGFDGM} & $\mathbf{L}= \|AH^{-1}A^T\|_1I$   & 2000/351k/993k
 & 5943.9 & 6122 &-&-& 2:52.9 \\
\cite{dualNewtonCG} & $\epsilon_i=10^{-4}, \mu=0.8,\sigma=0.3$   & 2000/351k/993k
 & 6464.1 & 20624 & 3961.9 & 12641 & 41:28.0 \\
\hline
Algorithm~\ref{alg:PGFDGM} & $\mathbf{L} = AH^{-1}A^T$  & 8000/1.41M/3.98M
 & -&-&2.0 & 2 & 9.4 \\
Algorithm~\ref{alg:DGFDGM} & $\mathbf{L}$ blk-diag comp. fr.  & 8000/1.41M/3.98M
 & 340.2 & 426&-&- & 44.6 \\
Algorithm~\ref{alg:DGFDGM} & $\mathbf{L}= \|AH^{-1}A^T\|_2I$   & 8000/1.41M/3.98M
 & 10583.4 & 10688&-&- & 17:05.3 \\
Algorithm~\ref{alg:DGFDGM} & $\mathbf{L}= \|AH^{-1}A^T\|_1I$   & 8000/1.41M/3.98M
 & 12801.2 & 12928&-&- & 20:40.2 \\
\end{tabular}
\label{tab:numEval}
\end{table*}

The proposed algorithm is evaluated by applying it to a randomly
generated systems
with a sparse dynamic interaction. The
dynamic interaction structure is decided using the method in
\cite[~\S 6.1]{Kraning} and the number of subsystems are 500, 2000,
and 8000 respectively. The resulting average degree of the generated
interconnections structures are 2.27, 2.23, and 2.23 respectively.
The number of states in each subsystem is randomly chosen from
the interval $\{10,11,\ldots,20\}$, the number of inputs are
three or four, and the control horizon is $N=10$. This
gives a total number of 87060, 350860, and 1405790 decision variables respectively. The
entries of the dynamics and input matrices are randomly
chosen from the intervals $[-0.7~1.3]$ and $[-1~1]$ respectively. Then
the dynamics matrix is re-scaled to get a spectral radius of
1.15. The states and inputs are upper and lower bounded by random
bounds generated from the intervals $[0.4~1]$ and $[-1~-0.4]$ respectively.
The state and input cost matrices are diagonal and each
diagonal entry is randomly chosen from the interval $[1~10^6]$.

The proposed algorithm is evaluated by comparing it to
fast dual decomposition, and the dual
Newton conjugate gradient (CG) method presented in \cite{dualNewtonCG}.
Fast dual
decomposition is achieved by setting $\mathbf{L}_i =
\|AH^{-1}A^T\|_2I$ for all $i$ in Algorithm~\ref{alg:DGFDGM}, where $A$
and $H$ are the global equality constraint 
and cost matrices respectively. This choice of $\mathbf{L}_i$ is
optimal if restricted to being a multiple of the identity matrix, and
if all $\mathbf{L}_i$ are restricted to be equal (as in fast dual
decomposition). However, this choice of $\mathbf{L}_i$ needs centralized computations,
which makes it unfair to call it a distributed method. We also compare
to fast dual gradient method using $\mathbf{L}_i =
\|AH^{-1}A^T\|_1I$ which satisfies $\|AH^{-1}A^T\|_1I\geq
\|AH^{-1}A^T\|_2I$. This choice of $\mathbf{L}_i$ can be computed
distributively with centralized coordination. We do not compare the
presented methods to standard dual decomposition with pure gradient steps,
since such methods are highly inferior. The dual Newton CG
method presented in \cite{dualNewtonCG} solves the dual problem using
a Newton method. The search direction is computed by solving the
resulting linear equations to some accuracy using distributed
conjugate gradient iterations. In each of these iterations, one local
and two global communications are performed. The Newton step-size is decided by a
distributed line search procedure that requires two global
communications for each function value comparison. In the algorithm,
the accuracy of the solution to the
linear system solved by the conjugate gradient method must be
specified. There is a trade-off between the number of iterations in the
CG-algorithm and the quality of the resulting search direction. If the
accuracy requirement is too low, we get close to a gradient direction,
which results in an expensive method that takes approximately gradient
steps. On the other hand, if the accuracy requirement is too high, too
many CG-iterations are performed in each iteration which gives a high
communication load. These algorithms are compared to the distributed
and parallel algorithms presented in this paper. For the parallel
algorithm, i.e.
Algorithm~\ref{alg:PGFDGM}, we choose
$\mathbf{L}=AH^{-1}A^T$ and pre-compute
the Cholesky factorization of this matrix for later use online. For
the distributed algorithm, i.e. Algorithm~\ref{alg:DGFDGM}, the
$\mathbf{L}_i$-matrices are computed based on
Algorithm~\ref{alg:distrInit}. In step 1) of Algorithm~\ref{alg:distrInit}, the
optimization problem \eqref{eq:findLocalL} is solved in each node $i$.
Finally, we note that all inner minimization problems in all
algorithms (also the Newton CG-algorithm) are solved using one max and
one min operation for each 
variable only. This is possible due to the diagonal structure of the
cost matrices and the since we have (soft) bound constraints only.

The evaluation in Table~\ref{tab:numEval} is obtained by generating
200 feasible random initial conditions from the state constraint set
for each of the systems. The corresponding optimal control problems
are solved using the
different algorithms, each utilizing 12 parallel cores. The first two
algorithms presented in Table~\ref{tab:numEval} for each problem batch
are the algorithms presented in this paper. The algorithms on row
three and four are fast dual decomposition with different step-sizes,
i.e. Algorithm~\ref{alg:DGFDGM} with $\mathbf{L}=\|AH^{-1}A^T\|_2I$ and
$\mathbf{L}=\|AH^{-1}A^T\|_1I$ respectively. The fifth and last row
for each problem batch contain results for the dual Newton CG method
in \cite{dualNewtonCG}. For each of these methods,
Table~\ref{tab:numEval} reports the average and
max number of local and global iterations, and the average execution
times for the 12 cores implementations. Due
to the very efficient implementation of the inner minimization
problem, the reported execution times are often dominated by the
execution time for the dual variable updates. In the general
situation with less efficient inner minimizations, the execution times
for the algorithms with many
local inner minimization problems would increase.

We start by comparing the two algorithms presented in this paper,
namely the distributed
Algorithm~\ref{alg:DGFDGM} and the parallel
Algorithm~\ref{alg:PGFDGM}. We first point out that
Algorithm~\ref{alg:DGFDGM} is fully distributed, both in
initialization and in execution, while Algorithm~\ref{alg:PGFDGM} is
initialized using centralized computations and is requires a global
communication structure. The number of communication rounds in
Algorithm~\ref{alg:PGFDGM} is substantially smaller than in
Algorithm~\ref{alg:DGFDGM}, but the communication in
Algorithm~\ref{alg:PGFDGM} is global. This is
due to the tighter quadratic upper bound
used in Algorithm~\ref{alg:PGFDGM}, i.e. $\mathbf{L}=AH^{-1}A^T$. Also, the average execution time
is smaller for Algorithm~\ref{alg:PGFDGM} in all examples for the 12 core
implementations. However, if using more computational units in the
algorithms, Algorithm~\ref{alg:DGFDGM} would outperform
Algorithm~\ref{alg:PGFDGM} in the ideal case where communication time
is neglected. This is due to the fully distributed structure of
Algorithm~\ref{alg:DGFDGM}. The possibility to achieve better execution times also in
practice using Algorithm~\ref{alg:DGFDGM}, hinges on the use of a very 
efficient synchronization and
communication protocol.

We also compare Algorithm~\ref{alg:DGFDGM} with block-diagonal
$\mathbf{L}$ as presented in this paper to fast dual decomposition with centralized
initialization, i.e. to Algorithm~\ref{alg:DGFDGM} with
$\mathbf{L}=\|AH^{-1}A^T\|_2I$, and to fast dual decomposition with
decentralized initialization, i.e. to Algorithm~\ref{alg:DGFDGM} with
$\mathbf{L}=\|AH^{-1}A^T\|_1I$. Table~\ref{tab:numEval} reveals that
the communicational burden
is greatly reduced using our algorithm. However, the complexity within each
iteration is slightly increased for 
Algorithm~\ref{alg:DGFDGM} with block-diagonal $\mathbf{L}$ compared
to fast dual decomposition. From the average execution times in Table~\ref{tab:numEval} we see
that this slight increase is by far compensated by the reduced number
of iterations. We also comment that if Table~\ref{tab:numEval} was
augmented with an entry for traditional dual decomposition, i.e. when
solving the dual problem using a standard gradient method, the
corresponding iteration count would be more that one order of
magnitude worse than for fast dual decomposition. This further
underlines the performance of our method.

Finally, we compare our algorithms to the recently proposed dual
Newton CG method in
\cite{dualNewtonCG}. As mentioned, the accuracy of the CG-gradient
method used to compute the search direction must be specified. We 
use $\epsilon_i=10^{-4}$ which is the least conservative accuracy
for which none of the
initial conditions gives too many Newton steps, where too many is in
the hundreds.
The main computational time as well as the main communicational burden
in the algorithm in \cite{dualNewtonCG} is spent on 
computing the search direction. This search direction is computed by
solving a linear system of equations using the conjugate gradient
method. This implies that an approximation of a system-wide inverse is
computed in every Newton iteration. For the problems considered here,
way too many CG-iterations are needed to compute a reasonable search
direction. This is revealed by 
Table~\ref{tab:numEval} that shows a significantly worse performance
of 
the method in \cite{dualNewtonCG} compared to our algorithms. For the
2000 sub-system problem, the average execution time was over 41
minutes, which gives a batch time of almost six days for all the 200
problems. The batch time for the 8000 sub-system problem would be in
the month range, which is why this is omitted from the comparison. The
performance evaluation is clear also without this table entry.

\section{Conclusions}

We have proposed a generalization of fast dual decomposition. In this
generalization, a quadratic upper bound to the
negative dual function with
different curvature in different directions is minimized in each step
in the algorithm. This
differs from traditional dual decomposition methods where the main
step is to minimize a quadratic upper bound to the negative dual
function that has the same curvature in all directions. This
generalization is made possible by the main contribution of this paper
that characterizes the set of matrices that can be used to describe
this quadratic upper bound. We propose on fully distributed algorithm
and one parallel algorithm, and we show that the distributed algorithm can be
initialized and reconfigured using distributed computations only. This
is traditionally not the case in dual decomposition where the norm of
a matrix that involve variables from all subsystems is used to compute
the optimal step size. The numerical evaluation shows that our
algorithms significantly outperform other distributed optimization
algorithms.

\bibliography{references}

\end{document}